\documentclass[twocolumn]{autart}    

\usepackage{graphicx}          

\usepackage{amsmath}
\usepackage{amssymb}

\newcommand*\diff{\mathop{}\!\mathrm{d}}%

\usepackage{graphics} 
\usepackage{epsfig} 
\usepackage{epstopdf}
\usepackage{subfigure}
\usepackage{color}

\graphicspath{{./Figures/}}
\allowdisplaybreaks[4]

\begin{document}

\begin{frontmatter}

\title{Robustness of constant-delay predictor feedback for in-domain stabilization of reaction-diffusion PDEs with time- and spatially-varying input delays\thanksref{footnoteinfo}} 

\thanks[footnoteinfo]{This work was supported by a research grant from Science Foundation Ireland (SFI) under grant number 16/RC/3872 and is co-funded under the European Regional Development Fund and by I-Form industry partners.\\
Corresponding author H.~Lhachemi.}

\author[CS]{Hugo Lhachemi}\ead{hugo.lhachemi@centralesupelec.fr},
\author[GIPSA-lab]{Christophe Prieur}\ead{christophe.prieur@gipsa-lab.fr},
\author[ICL]{Robert Shorten}\ead{r.shorten@imperial.ac.uk},               

\address[CS]{L2S, CentraleSup{\'e}lec, 91192 Gif-sur-Yvette,
France}  
\address[GIPSA-lab]{Universit{\'e} Grenoble Alpes, CNRS, Grenoble-INP, GIPSA-lab, F-38000, Grenoble, France}             
\address[ICL]{Dyson School of Design Engineering, Imperial College London, London, U.K}        

\begin{keyword}                           
Delayed distributed actuation, Spatially-varying delay, Distributed parameter systems, Predictor feedback, Reaction-diffusion equation              
\end{keyword}                             

\begin{abstract}                          
This paper discusses the in-domain feedback stabilization of reaction-diffusion PDEs with Robin boundary conditions in the presence of an uncertain time- and spatially-varying delay in the distributed actuation. The proposed control design strategy consists of a constant-delay predictor feedback designed based on the known nominal value of the control input delay and is synthesized on a finite-dimensional truncated model capturing the unstable modes of the original infinite-dimensional system. By using a small-gain argument, we show that the resulting closed-loop system is exponentially stable provided that the variations of the delay around its nominal value are small enough. The proposed proof actually applies to any distributed-parameter system associated with an unbounded operator that 1) generates a $C_0$-semigroup on a weighted space of square integrable functions over a compact interval; and 2) is self-adjoint with compact resolvent.
\end{abstract}

\end{frontmatter}

\section{Introduction}
Stabilization of open-loop unstable partial differential equations (PDEs) in the presence of delays has attracted much attention in the recent years. A first class of problems deals with the feedback stabilization of PDEs in the presence of a state-delay~\cite{fridman2009exponential,hashimoto2016stabilization,kang2017boundary,kang2017boundaryIFAC,kang2018boundary,lhachemi2019boundary,lhachemi2019boundaryISS,solomon2015stability}. In this paper, we are concerned with a second class of problem, namely: the feedback stabilization of PDEs in the presence of a delay in the control input~\cite{guzman2019stabilization,krstic2009control,lhachemi2019feedback,lhachemi2019lmi,lhachemi2020exponential,lhachemi2019control,lhachemi2019pi,nicaise2008stabilization,nicaise2007stabilization,nicaise2009stability,prieur2018feedback,qi2019stabilization}. One of the very first contributions in this field was reported in~\cite{krstic2009control}. In this work, the problem of boundary feedback stabilization of an unstable reaction-diffusion equation under a constant input delay was tackled via a backstepping transformation. More recently, the same problem was investigated in~\cite{prieur2018feedback} by adopting a different control design approach. Inspired by the early work~\cite{russell1978controllability} and the later developments reported in~\cite{coron2004global,coron2006global}, the authors synthesized a predictor feedback on a finite-dimensional model capturing the unstable modes of the original infinite-dimensional system. The stability property of the resulting closed-loop infinite-dimensional system was obtained via the study of a Lyapunov function. It was shown in~\cite{guzman2019stabilization} that this approach is not limited to reaction-diffusion systems but can also be applied to the boundary feedback stabilization of a linear Kuramoto-Sivashinsky equation under a constant input delay. This approach was generalized to the boundary stabilization of a class of diagonal infinite-dimensional systems in~\cite{lhachemi2019feedback,lhachemi2019control} for constant input delays and then in~\cite{lhachemi2019lmi,lhachemi2020exponential} for fast time-varying input delays.

Most of the approaches reported in the literature deal with boundary control inputs only. Very few reported works are concerned with the in-domain stabilization of PDEs in the presence of a long delay in the control input. In this domain, the recent work~\cite{qi2019stabilization} tackles the in-domain stabilization of an unstable reaction-diffusion equation with Dirichlet boundary conditions and a constant delay in the in-domain control input. The reported control design strategy takes advantage of a backstepping transformation and involves technical challenges in the stability analysis due to the occurrence of kernel functions presenting singularities. 

From a practical perspective, it is worth noting that input delays are generally uncertain and possibly time-varying. In this context, the study of the robustness of the proposed control strategies with respect to delay mismatches is of paramount importance. The case of a distributed actuation scheme is even more complex since spatially-varying delays can arise due to network and transport effects that may vary among different spatial regions. A first example of this situation occurs in the context of biological systems and population dynamics~\cite{schley1999linear}. In such a situation, delays effects are ubiquitous due to reaction or maturation times induced either by natural processes or exogenous inputs acting in a feedback loop. An example of the latter can be found in the context of epidemic dynamics~\cite{wang2020regional} in which control inputs take the form of either medical prescriptions (medicines, vaccination), social distancing measures or physical restrictions (confinement, partial limitation of people fluxes). In this setting, spatially-varying delays appear in the application of the measures due to the combination of incubation periods~\cite{Guan2020.07.27.20161430} and specific regional characteristics. A second example occurs in the context of thermonuclear fusion with Tokamaks~\cite{mavkov2017distributed}. The objective of these devices is to control the plasma in a torus in order to, ultimately, achieve controlled thermonuclear fusion. In this setting, one of the control design objectives is to regulate the temperature of the plasma's electrons described by a diffusion equation. The distributed control input takes the form of the total electron heating power density and is actually implemented by a set of neutral-beam injection and radio frequency antennas; see in particular~\cite[Eq. (2)]{mavkov2017distributed} that is a diffusion equation with one distributed control input. It is reported in~\cite{argomedo2014safety} that, due to network effects, delays of around 100~ms are introduced in the feedback loop. Uncertain time- and spatially-varying delay occur due to network effects and the multiplicity of the devices used to generate the heating power density. A last example occurs in the context of the stabilization of fronts in a reaction-diffusion system with possible application to chemical reactors~\cite{smagina2002stabilization}. However, to the best of our knowledge, the design and/or robustness analysis of control strategies with respect to possibly spatially-varying delays is still an open problem. The present study is a first step into that research direction.

This paper is concerned with the feedback stabilization of an unstable reaction-diffusion equation with Robin boundary conditions in the presence of an uncertain time- and spatially-varying delay in the distributed control input. Motivated by~\cite{prieur2018feedback}, the proposed control strategy relies on a constant-delay predictor feedback synthesized on a finite-dimensional truncated model capturing the unstable modes of the original infinite-dimensional system. In essence, this approach is similar to the one reported in~\cite{lhachemi2019lmi} with application to the boundary control of a class of diagonal abstract boundary control systems. However, we point out that the spatially-varying nature of the delay in the control input brings new challenges that do not allow the replication of the proof of stability reported in~\cite{lhachemi2019lmi}. This is because while time and space variables where fully uncoupled in~\cite{lhachemi2019lmi}, the spatially-varying nature of the delay considered in this present work introduces a strong coupling between time and space variables. Consequently, a dedicated stability analysis is required. Inspired by the early work~\cite{karafyllis2013delay} dealing with the robustness of constant-delay predictor feedback w.r.t. uncertain and time-varying input delays for finite-dimensional systems (see also~\cite{lhachemi2020exponential} in the context of input-to-state stabilization), this analysis is carried out in this paper via a small gain argument. We show that the constant-delay predictor feedback achieves the exponential stabilization of the closed-loop infinite-dimensional system provided that the deviations of the uncertain time- and spatially-varying delay around its nominal value are small enough. The derived proof applies to any distributed parameter system associated with an unbounded operator that 1) generates a $C_0$-semigroup on a weighted space of square integrable functions on a compact interval; and 2) is self-adjoint with compact resolvent. This includes, e.g., the linear Kuramoto-Sivashinsky equation studied in~\cite{guzman2019stabilization}.

The remainder of this paper is organized as follows. The problem setting and the proposed control strategy are reported in Section~\ref{sec: problem setting}. Then, the stability analysis is carried out in Section~\ref{sec: stability analysis}. The numerical illustration of the obtained results is reported in Section~\ref{sec: numerical application}. Finally, concluding remarks are provided in Section~\ref{sec: conclusion}.

\textbf{Notation.} The sets of non-negative integers, real, and non-negative real numbers are denoted by $\mathbb{N}$, $\mathbb{R}$, and $\mathbb{R}_+$, respectively. The set of $n$-dimensional vectors over $\mathbb{R}$ is denoted by $\mathbb{R}^n$ and is endowed with the Euclidean norm $\Vert x \Vert = \sqrt{x^* x}$. The set of $n \times m$ matrices over $\mathbb{R}$ is denoted by $\mathbb{R}^{n \times m}$ and is endowed with the induced norm denoted by $\Vert\cdot\Vert$. For any $t_0>0$, we say that $\varphi \in \mathcal{C}^0(\mathbb{R};\mathbb{R})$ is a \textit{transition signal over $[0,t_0]$} if $0 \leq \varphi \leq 1$, $\left. \varphi \right\vert_{(-\infty,0]} = 0$, and $\left. \varphi \right\vert_{[t_0,+\infty)} = 1$.

\section{Problem setting and control design strategy}\label{sec: problem setting}

\subsection{Problem setting}

\subsubsection{Abstract system}\label{subsec: abstract system}

We consider the real state-space $\mathcal{H} = L_\rho^2(0,1)$ for some $0 < \rho \in \mathcal{C}^0([0,1];\mathbb{R})$, i.e. the space of square integrable functions over $(0,1)$ endowed with the weighted\footnote{The introduction of the weighting function $\rho$ is motivated by the study of the reaction-diffusion equation described in Subsection~\ref{subsec: example reaction-diffusion}.} inner product $\left< f , g\right> = \int_0^1 \rho(\xi) f(\xi) g(\xi) \diff\xi$. The associated norm is denoted by $\Vert \cdot \Vert_\mathcal{H}$. We recall that this structure defines a separable real Hilbert space. Let $\mathcal{A} : \mathcal{D}(\mathcal{A}) \subset \mathcal{H} \rightarrow \mathcal{H}$ be the generator of a $C_0$-semigroup $T(t)$. We further assume that $\mathcal{A}$ is self-adjoint with compact resolvent. In this context the following result is standard, see e.g. \cite[Chap.~6]{brezis2010functional} and \cite[Sec~A.4.2]{Curtain2012}. The eigenvalues $(\lambda_n)_{n \geq 1}$ of $\mathcal{A}$ are all real with finite multiplicity, can be sorted such that they form a non-increasing sequence with $\lambda_n \rightarrow - \infty$ when $n \rightarrow + \infty$, and the associated eigenvectors $(e_n)_{n \geq 1}$ can be selected to form a Hilbert basis of $\mathcal{H}$.

Our starting point is the abstract system:
\begin{subequations}\label{eq: abstract form}
\begin{align}
\dfrac{\mathrm{d}X}{\mathrm{d}t}(t) & = \mathcal{A}X(t) + v(t) \label{eq: abstract form - DE} \\
X(0) & = X_0 \label{eq: abstract form - IC}
\end{align}
\end{subequations}
for $t>0$. Here $X(t) \in \mathcal{H}$ is the state-vector and $X_0 \in\mathcal{H}$ is the initial condition. We assume that the distributed feedback control $u(t) \in \mathcal{H}$ is related to $v(t) \in\mathcal{H}$ by
\begin{equation*}
[v(t)](\xi) = [u(t-D(t,\xi))](\xi)
\end{equation*}
with $D \in \mathcal{C}^0(\mathbb{R}_+ \times [0,1] ; \mathbb{R})$ a time- and spatially-varying delay that satisfies $\vert D - D_0 \vert \leq \delta$ where $D_0 > 0$ and $\delta \in (0,D_0)$ are known given constants. Constant $D_0 > 0$ is referred to as the nominal value of the delay $D$ while $\delta > 0$ stands for its maximal amplitude of variation around $D_0$. The system is assumed uncontrolled for negative times, i.e., $[u(t)](\xi) = 0$ for $t < 0$ and $\xi\in(0,1)$. The objective is to design the feedback control $u$, taking the form of a state-feedback of the system trajectory $X$, such that the closed-loop system is exponentially stable.

\subsubsection{Example 1: reaction-diffusion equation}\label{subsec: example reaction-diffusion}

The abstract formulation as previously described is motivated by the study of the in-domain feedback stabilization of the following reaction-diffusion equation with Robin boundary conditions:
\begin{subequations}\label{eq: pb setting}
\begin{align}
& y_t(t,\xi) = \dfrac{1}{\rho(\xi)} (p y_{\xi})_\xi (t,\xi) + \dfrac{q(\xi)}{\rho(\xi)} y(t,\xi)\label{eq: pb setting - PDE} \\
& \phantom{y_t(t,\xi) =}\, + u(t-D(t,\xi),\xi) \nonumber \\
& \cos(\theta_1) y(t,0) - \sin(\theta_1) y_\xi(t,0) = 0 \label{eq: pb setting - BC1} \\
& \cos(\theta_2) y(t,1) + \sin(\theta_2) y_\xi(t,1) = 0 \label{eq: pb setting - BC2} \\
& y(0,\xi) = y_0(\xi), \label{eq: pb setting - IC}
\end{align}
\end{subequations}
for $t > 0$ and $\xi \in (0,1)$. Here we have $\rho,q \in\mathcal{C}^0([0,1];\mathbb{R})$, $p \in\mathcal{C}^1([0,1];\mathbb{R})$, $\rho,p>0$, and $\theta_1,\theta_2 \in [0,2\pi)$. In this setting, $u : [-D_0-\delta,+\infty) \times (0,1) \rightarrow \mathbb{R}$, with $u(t,\cdot) = 0$ for $t < 0$, is the in-domain control input. This input is subject to the uncertain time- and spatially-varying continuous input delay $D : \mathbb{R}_+ \times [0,1] \rightarrow \mathbb{R}$ with $\vert D - D_0 \vert \leq \delta$ where $D_0 > 0$ and $\delta \in (0,D_0)$ are given constants. Finally, $y_0 : (0,1) \rightarrow \mathbb{R}$ stands for the initial condition.

The reaction-diffusion system (\ref{eq: pb setting}) can be written in the abstract form (\ref{eq: abstract form}) by using the real state-space $\mathcal{H} = L_\rho^2(0,1)$. In this case, we have the operator $\mathcal{A}f = \frac{1}{\rho} (p f')' + \frac{q}{\rho} f \in\mathcal{H}$ defined on the domain $\mathcal{D}(\mathcal{A}) = \{ f \in H^2(0,1) \,:\, \cos(\theta_1) f(0) - \sin(\theta_1) f'(0) = 0 ,\, \cos(\theta_2) f(1) + \sin(\theta_2) f'(1) = 0 \}$, the state-vector $X(t) = y(t,\cdot) \in \mathcal{H}$, the distributed function $v(t) = u(t-D(t,\cdot),\cdot) \in\mathcal{H}$ with control input $u(t,\cdot) \in\mathcal{H}$, and the initial condition $X_0=y_0 \in \mathcal{H}$. Recalling that $\mathcal{A}$ generates a $C_0$-semigroup $T(t)$ on $\mathcal{H}$ and that $\mathcal{A}$ is self-adjoint with compact resolvent (see, e.g., \cite[Sec.~8.6]{renardy2006introduction} and~\cite{delattre2003sturm}), the context of the abstract form (\ref{eq: abstract form}) applies to this system.

\begin{rem}
The stabilization of (\ref{eq: pb setting}) in the case of constant functions $\rho,p,q$, a constant and known delay $D$, and for Dirichlet boundary conditions ($\theta_1=\theta_2=0$), has been investigated in~\cite{qi2019stabilization} via a backstepping design.
\end{rem}

\subsubsection{Example 2: linear Kuramoto-Sivashinsky equation}

An other example of a PDE system fitting within the abstract form (\ref{eq: abstract form}) is the linear Kuramoto-Sivashinsky equation studied in~\cite{guzman2019stabilization}:
\begin{subequations}\label{eq: linear Kuramoto-Sivashinsky equation}
\begin{align}
& y_t(t,\xi) + y_{\xi\xi\xi\xi}(t,\xi) + \lambda y_{\xi\xi}(t,\xi) = u(t-D(t,\xi),\xi) \\
& y(t,0) = y(t,1) = y_\xi(t,0) = y_\xi(t,1) = 0 \\
& y(0,\xi) = y_0(\xi) ,
\end{align}
\end{subequations}
for $t > 0$ and $\xi \in (0,1)$. Here we have $\lambda > 0$. As in the previous setting, $u$ is the in-domain control input, $D$ is a time- and spatially-varying delay, and $y_0$ is the initial condition.

The linear Kuramoto-Sivashinsky equation (\ref{eq: linear Kuramoto-Sivashinsky equation}) can be written as (\ref{eq: abstract form}) by introducing the real state-space $\mathcal{H} = L^2(0,1)$, the operator $\mathcal{A}f = - f'''' - \lambda f'' \in\mathcal{H}$ defined on the domain $\mathcal{D}(\mathcal{A}) = H^4(0,1) \cap H_0^2(0,1)$, the state-vector $X(t) = y(t,\cdot) \in \mathcal{H}$, the distributed function $v(t) = u(t-D(t,\cdot),\cdot) \in\mathcal{H}$ with control input $u(t,\cdot) \in\mathcal{H}$, and the initial condition $X_0=y_0 \in \mathcal{H}$. The fact that $\mathcal{A}$ is self-adjoint, has compact resolvent, and generates a $C_0$-semigroup, is reported, e.g., in~\cite{cerpa2017control}.

\subsection{Control design strategy}
Assuming that the control input $u$ is such that\footnote{This regularity will be assessed in Subsection~\ref{subsec: well-posedness} based on the forthcoming control strategy.} $v \in \mathcal{C}^0(\mathbb{R}_+;\mathcal{H})$, the mild solution $X\in\mathcal{C}^0(\mathbb{R}_+;\mathcal{H})$ of (\ref{eq: abstract form}) is uniquely defined by~\cite[Def.~3.1.4 and Lem.~3.1.5]{Curtain2012}
\begin{equation}\label{eq: def mild solution}
X(t) = T(t)X_0 + \int_0^t T(t-s) v(s) \,\mathrm{d}s .
\end{equation}
We introduce $x_n(t) = \left< X(t) , e_n \right>$ the coefficients of projection of $X(t)$ onto the Hilbert basis $(e_n)_{n \geq 1}$. Then we have $X(t) = \sum_{n \geq 1} x_n(t) e_n$ and $\Vert X(t) \Vert_\mathcal{H}^2 = \sum_{n \geq 1} \vert x_n(t) \vert^2$ for all $t \geq 0$. Since $\mathcal{A}e_n = \lambda_n e_n$, we have that $T(t) e_n = e^{\lambda_n t} e_n$. Thus, we obtain from (\ref{eq: def mild solution}) that
\begin{equation*}
x_n(t) = e^{\lambda_n t} x_n(0) + \int_0^t e^{\lambda_n (t-s)} \left< v(s) , e_n \right> \,\mathrm{d}s .
\end{equation*}
As $v$ is continuous, this shows that $x_n \in \mathcal{C}^1(\mathbb{R}_+;\mathbb{R})$ and satisfies the ODE
\begin{equation*}
\dot{x}_n(t) = \lambda_n x_n(t) + \left< v(t) , e_n \right>
\end{equation*}
for all $t \geq 0$. Considering $D_0 > 0$ a nominal value of the delay $D$ as described in Subsection~\ref{subsec: abstract system}, we define a nominal delayed control input $v_0(t) = u(t-D_0)$. We also introduce the coefficients of projection $v_n(t) = \left< v(t) , e_n \right>$ and $v_{0,n}(t) = \left< v_0(t) , e_n \right>$, and the residual term $\Delta_n(t) = v_n(t) - v_{0,n}(t) = \left< v(t) - v_{0}(t) , e_n \right>$. Then we have
\begin{equation}\label{eq: spectral reduction - scalar ODE}
\dot{x}_n(t) = \lambda_n x_n(t) + v_{0,n}(t) + \Delta_n(t)
\end{equation}
for all $t \geq 0$.

Let $N \geq 1$ and $\gamma > 0$ be such that $\lambda_n \leq - \gamma$ for all $n \geq N+1$. We consider the following structure for the control input:
\begin{equation}\label{eq: structure control input u}
[u(t)](\xi) = \sum\limits_{k=1}^{N} w_k(t) e_k(\xi)
\end{equation}
with $w_k(t) \in \mathbb{R}$ to be defined.  In particular, we have
\begin{subequations}\label{eq: structure control input v and v0}
\begin{align}
[v(t)](\xi) & = \sum\limits_{k=1}^{N} w_k(t-D(t,\xi)) e_k(\xi) , \label{eq: structure control input v} \\
[v_0(t)](\xi) & = \sum\limits_{k=1}^{N} w_k(t-D_0) e_k(\xi) . \label{eq: structure control input v0}
\end{align}
\end{subequations}
Then we obtain from (\ref{eq: spectral reduction - scalar ODE}) that
\begin{equation}\label{eq: spectral reduction - scalar ODE - 1}
\dot{x}_n(t) = \lambda_n x_n(t) + w_n(t-D_0) + \Delta_n(t)
\end{equation}
for $1 \leq n \leq N$, while
\begin{equation}\label{eq: spectral reduction - scalar ODE - 2}
\dot{x}_n(t) = \lambda_n x_n(t) + v_n(t)
\end{equation}
for $n \geq N+1$.

\begin{rem}
As it can be seen from (\ref{eq: structure control input v}), the spatially-varying nature of the input delay introduces a strong coupling between the time and space variables. A decoupling is obtained only in the case of a delay that is uniform throughout the spatial domain, i.e., $D(t,\xi) = D_u(t)$. In that case, (\ref{eq: structure control input v}) reduces to $[v(t)](\xi) = \sum_{k=1}^{N} w_k(t-D_u(t)) e_k(\xi)$. This implies the following simplifications: $v_n(t) = w_n(t-D_u(t))$ and $\Delta_n (t) = w_n(t-D_u(t)) - w_n(t-D_0)$ for $n \leq N$ while $v_n(t) = 0$ for $n \geq N+1$.
\end{rem}

Introducing
\begin{align*}
x(t) & = \begin{bmatrix} x_1(t) & \ldots & x_{N}(t) \end{bmatrix}^\top \in \mathbb{R}^{N} , \\
w(t) & = \begin{bmatrix} w_1(t) & \ldots & w_{N}(t) \end{bmatrix}^\top \in \mathbb{R}^{N} , \\
\Delta(t) & = \begin{bmatrix} \Delta_1(t) & \ldots & \Delta_{N}(t) \end{bmatrix}^\top \in \mathbb{R}^{N} , \\
\Lambda & = \mathrm{diag}(\lambda_1,\ldots,\lambda_N) \in \mathbb{R}^{N \times N} ,
\end{align*}
we obtain that
\begin{equation}\label{eq: spectral reduction - vect ODE}
\dot{x}(t) = \Lambda x(t) + w(t-D_0) + \Delta(t)
\end{equation}
for all $t \geq 0$. From (\ref{eq: structure control input u}), we have $\Vert u(t) \Vert_\mathcal{H} = \Vert w(t) \Vert$.

The control design strategy consists of the design of a constant-delay predictor feedback in the nominal configuration $D(t,\xi) = D_0$ for which (\ref{eq: spectral reduction - vect ODE}) reduces to $\dot{x}(t) = \Lambda x(t) + w(t-D_0)$. Thus, the control scheme takes the form of the classical constant-delay predictor feedback:
\begin{equation}\label{eq: control input w}
w(t) = \varphi(t) K \left\{ x(t) + \int_{t-D_0}^{t} e^{(t-D_0-s)\Lambda} w(s) \,\mathrm{d}s \right\} ,
\end{equation}
where $K \in \mathbb{R}^{N \times N}$ is a feedback gain such that $A_\mathrm{cl} = \Lambda + e^{-D_0 \Lambda} K$ is Hurwitz and $\varphi \in \mathcal{C}^0(\mathbb{R};\mathbb{R})$ is a transition signal\footnote{See notation section.} over $[0,t_0]$ for some arbitrarily given $t_0 > 0$. In particular, we have $w(t) = 0$ and hence  $u(t) = 0$ for $t \leq 0$. The existence and uniqueness of a function $w$ that is solution of the implicit equation (\ref{eq: control input w}) has been investigated in~\cite{bresch2018new}. See the proof of Lemma~\ref{lem: well-posedness} for details.

The objective of the remainder of this paper is to show the following robustness result: the constant-delay predictor feedback (\ref{eq: control input w}) achieves the exponential stabilization of (\ref{eq: abstract form}) with command input (\ref{eq: structure control input u}) for small enough deviations of the time- and spatially-varying delay $D(t,\xi)$ around its nominal value $D_0$.

\begin{rem}
For a given desired closed-loop matrix $A_\mathrm{cl} \in \mathbb{R}^{N \times N}$, the corresponding feedback gain $K \in \mathbb{R}^{N \times N}$ is given by $K = e^{D_0 \Lambda} (A_\mathrm{cl}-\Lambda)$.
\end{rem}

\begin{rem}
The transition signal $\varphi$ appearing in (\ref{eq: control input w}) is used to ensure a continuous transition from open-loop ($t < 0$) to closed-loop ($t \geq 0$). In particular, recalling that $[v(t)](\xi) = [u(t-D(t,\xi))](\xi)$, this transition signal prevents the occurrence of jumps in the distributed signal $v(t)$ at times $t \geq 0$ for which the function $t \mapsto t- D(t,\xi)$ crosses $0$ while avoiding the introduction of compatibility conditions restricting the set of admissible initial conditions $X_0 \in \mathcal{H}$. This continuous behavior will be used in the well-posedness assessment; see the proof of Lemma~\ref{lem: well-posedness} for details.
\end{rem}

\subsection{Statement of the main result}

The main result of this paper is stated below.

\begin{thm}\label{thm: main result}
Let the real state-space $\mathcal{H} = L_\rho^2(0,1)$ for some $0 < \rho \in \mathcal{C}^0([0,1];\mathbb{R})$. Let $\mathcal{A} : \mathcal{D}(\mathcal{A}) \subset \mathcal{H} \rightarrow \mathcal{H}$ be a self-adjoint operator with compact resolvent and which is the generator of a $C_0$-semigroup. Let an integer $N \geq 1$ be such that $\lambda_{N+1} < 0$. Let $D_0 > 0$ be a given nominal delay. Let $K \in \mathbb{R}^{N \times N}$ be a feedback gain such that $A_\mathrm{cl} = \Lambda + e^{-D_0 \Lambda} K$ is Hurwitz. Let $M \geq 1$ and $\sigma > 0 $ be such that $\Vert e^{A_\mathrm{cl} t} \Vert \leq M e^{-\sigma t}$ for all $t \geq 0$. We denote by $K_k$ the $k$-th line of the feedback gain $K$. Let $\delta \in (0,D_0)$ be such that
\begin{equation}\label{eq: small gain condition}
\dfrac{\max( M , e^{\Vert A_\mathrm{cl} \Vert \delta} ) \sqrt{N}}{\sigma} \sum\limits_{k=1}^{N} \Vert K_k \Vert \left\{ (e^{\Vert A_\mathrm{cl} \Vert \delta} - 1) + \sigma\delta e^{\sigma\delta} \right\} < 1 .
\end{equation}
Let $\varphi \in \mathcal{C}^0(\mathbb{R};\mathbb{R})$ be a transition signal over $[0,t_0]$ for some given $t_0 > 0$. Then there exist constants $\kappa,C>0$ such that, for any initial condition $X_0 \in \mathcal{H}$ and any delay $D \in \mathcal{C}^0(\mathbb{R}_+\times[0,1];\mathbb{R})$ with $\vert D - D_0 \vert \leq \delta$, the mild solution $X\in\mathcal{C}^0(\mathbb{R}_+;\mathcal{H})$ of the closed-loop system composed of (\ref{eq: abstract form}), (\ref{eq: structure control input u}), and (\ref{eq: control input w}) satisfies
\begin{equation*}
\Vert X(t) \Vert_\mathcal{H} + \Vert u(t) \Vert_\mathcal{H} \leq C e^{-\kappa t} \Vert X_0 \Vert_\mathcal{H}
\end{equation*}
for all $t \geq 0$.
\end{thm}

\begin{rem}
As the left-hand side of (\ref{eq: small gain condition}) is equal to zero when evaluated at $\delta = 0$, the existence of $\delta > 0$ such that (\ref{eq: small gain condition}) holds is ensured by a continuity argument. This shows that the constant-delay predictor feedback synthesized based on the nominal value $D_0$ of the time- and spatially-varying delay $D(t,\xi)$ ensures the exponential stability of the resulting closed-loop system for delays with deviations $\delta$ around the nominal value $D_0$ that are small enough. In this context, (\ref{eq: small gain condition}) stands for an explicit sufficient condition on the admissible values of $\delta > 0$.
\end{rem}

\begin{rem}
The first part of the proof of Thm.~\ref{thm: main result} consists of the study of the robustness of the constant-delay predictor feedback with respect to delay mismatches in the context of the finite-dimensional system (\ref{eq: spectral reduction - vect ODE}). Note that similar problems were investigated in~\cite{bekiaris2013robustness,karafyllis2013delay,krstic2008lyapunov,lhachemi2019lmi,li2014robustness,selivanov2016predictor} either in the case of constant or time-varying input delays. However, due to the spatially varying nature of the delay considered in this work, the above results do not apply because of the occurrence of the $\Delta(t)$ term in (\ref{eq: spectral reduction - vect ODE}). Hence, a dedicated stability analysis, taking into account the spatially-varying nature of the delay, is required.
\end{rem}

\begin{rem}
The stability result stated by Theorem~\ref{thm: main result} holds in $L^2$-norm. In the particular case of the reaction-diffusion equation
\begin{equation*}
y_t(t,\xi) = y_{\xi\xi}(t,\xi) + c(\xi) y(t,\xi) + u(t-D(t,\xi),\xi)
\end{equation*}
with Dirichlet boundary conditions $y(t,0) = y(t,1) = 0$ and where $c \in L^\infty(0,1)$, the classical solutions (obtained, e.g., for $y_0 \in H^2(0,1) \cap H_0^1(0,1)$, $D \in \mathcal{C}^1(\mathbb{R}_+\times[0,1];\mathbb{R})$, and $\varphi \in \mathcal{C}^1(\mathbb{R};\mathbb{R})$) of the closed-loop system are exponentially stable in $H^1$-norm. This result essentially relies on the identity:
\begin{equation*}
\Vert f \Vert_{H_0^1(0,1)}^2
= \int_0^1 c(\xi) f(\xi)^2 \diff\xi
- \sum\limits_{n \geq 1} \lambda_n \langle f , e_n \rangle_{L^2(0,1)}^2
\end{equation*}
which holds for any $f \in H^2(0,1) \cap H_0^1(0,1)$; see \cite[Eq.~42]{prieur2018feedback} for a detailed proof. Based on the stability result of Theorem~\ref{thm: main result} and using a similar approach to the one reported in Subsection~\ref{subsec: stab residual dynamics} to estimate $\sum_{n \geq N+1} \vert \lambda_n \vert \vert x_n(t) \vert^2 $, the claimed stability estimate in $H^1$-norm follows.
\end{rem}

\section{Proof of the main result}\label{sec: stability analysis}

This section is devoted to the proof of the main result of this paper, namely: Theorem~\ref{thm: main result}.

\subsection{Well-posedness}\label{subsec: well-posedness}

We first assess the well-posedness of the closed-loop system dynamics.

\begin{lem}\label{lem: well-posedness}
For any initial condition $X_0 \in \mathcal{H}$ and any delay $D \in \mathcal{C}^0(\mathbb{R}_+\times[0,1];\mathbb{R})$ with $\vert D - D_0 \vert \leq \delta < D_0$, there exists a unique mild solution $X\in\mathcal{C}^0(\mathbb{R}_+;\mathcal{H})$ of the closed-loop system composed of (\ref{eq: abstract form}), (\ref{eq: structure control input u}), and (\ref{eq: control input w}). Moreover, the control input satisfies $u \in \mathcal{C}^0([-D_0-\delta,+\infty);\mathcal{H})$ as well as $v \in \mathcal{C}^0(\mathbb{R}_+;\mathcal{H})$ with $w \in \mathcal{C}^0([-D_0-\delta,+\infty);\mathbb{R}^N)$.
\end{lem}

\textbf{Proof.}
Let $X_0 \in \mathcal{H}$ and $D \in \mathcal{C}^0(\mathbb{R}_+ \times[0,1];\mathbb{R})$ with $\vert D - D_0 \vert \leq \delta < D_0$. We show by induction that, for any $k \geq 1$, the mild solution $X\in\mathcal{C}^0([0,k(D_0-\delta)];\mathcal{H})$ given by (\ref{eq: def mild solution}) is well and uniquely defined with $u \in \mathcal{C}^0([-D_0-\delta,k(D_0-\delta)];\mathcal{H})$ and $v \in \mathcal{C}^0([0,k(D_0-\delta)];\mathcal{H})$ where $w \in \mathcal{C}^0([-D_0-\delta,k(D_0-\delta)];\mathbb{R}^N)$ is the unique solution of (\ref{eq: control input w}) over the time interval $[-D_0-\delta,k(D_0-\delta)]$.

\emph{Initialization.} For $0 \leq t \leq D_0 - \delta$, we have that $t - D(t,\xi) \leq t - (D_0 - \delta) \leq 0$ hence $v(t)=0$. Then we have $X(t) = T(t) X_0$ for all $0 \leq t \leq D_0 - \delta$, yielding $X\in\mathcal{C}^0([0,D_0-\delta];\mathcal{H})$. In particular $x\in\mathcal{C}^0([0,D_0-\delta];\mathbb{R}^N)$ and the control input $w$ solution of the fixed-point equation (\ref{eq: control input w}) is well and uniquely defined (see~\cite{bresch2018new} for details), and we have $w \in \mathcal{C}^0([-D_0-\delta,D_0-\delta];\mathbb{R}^N)$. Finally, we infer from (\ref{eq: structure control input u}) that $u \in \mathcal{C}^0([-D_0-\delta,D_0-\delta];\mathcal{H})$.

\emph{Induction.} Assume that the property holds true for a given integer $k \geq 1$. For $0 \leq t \leq (k+1)(D_0 - \delta)$, we have that $t - D(t,\xi) \leq t - (D_0 - \delta) \leq k(D_0 - \delta)$. Thus $v$ over the time interval $[0,(k+1)(D_0-\delta)]$ only depends on the known control input $u$ for times in the interval $[-D_0-\delta,k(D_0-\delta)]$. We need to show that $v \in \mathcal{C}^0([0,(k+1)(D_0-\delta)];\mathcal{H})$. First, as $w_k$ is continuous on $[-D_0-\delta,k(D_0-\delta)]$ and $D$ is continuous on $\mathbb{R}_+ \times [0,1]$, we obtain that $w_k(t-D(t,\cdot)) \in L^\infty(0,1)$ for any $t \in [0,(k+1)(D_0-\delta)]$. Then we obtain from (\ref{eq: structure control input v}) that $v(t) \in L_\rho^2(0,1)$ for any $t \in [0,(k+1)(D_0-\delta)]$. Now we note from (\ref{eq: structure control input v}) that, for any $\tau,t \in [0,(k+1)(D_0-\delta)]$,
\begin{align*}
& \Vert v(\tau) - v(t) \Vert_\mathcal{H} \\
& \leq \sum\limits_{k=1}^{N} \Vert \{ w_k(\tau-D(\tau,\cdot)) - w_k(t-D(t,\cdot)) \} e_k \Vert_\mathcal{H}
\end{align*}
with
\begin{align*}
& \Vert \{ w_k(\tau-D(\tau,\cdot)) - w_k(t-D(t,\cdot)) \} e_k \Vert_\mathcal{H}^2 \\
& = \int_0^1 \rho(\xi) \vert w_k(\tau-D(\tau,\xi)) - w_k(t-D(t,\xi)) \vert^2 e_k(\xi)^2 \,\mathrm{d}\xi \\
& \underset{\tau \rightarrow t}{\longrightarrow} 0
\end{align*}
by the Lebesgue dominated convergence theorem~\cite{folland2013real}. We have shown that $v \in \mathcal{C}^0([0,(k+1)(D_0-\delta)];\mathcal{H})$. Thus, using (\ref{eq: def mild solution}), the mild solution $X\in\mathcal{C}^0([0,k(D_0-\delta)];\mathcal{H})$ is uniquely extended as a function $X\in\mathcal{C}^0([0,(k+1)(D_0-\delta)];\mathcal{H})$. In particular $x\in\mathcal{C}^0([0,(k+1)(D_0-\delta)];\mathbb{R}^N)$ and the control input $w$ solution of the fixed-point equation (\ref{eq: control input w}) is well and uniquely defined (see~\cite{bresch2018new} for details), and we have $w \in \mathcal{C}^0([-D_0-\delta,(k+1)(D_0-\delta)];\mathbb{R}^N)$. Finally, we infer from (\ref{eq: structure control input u}) that $u \in \mathcal{C}^0([-D_0-\delta,(k+1)(D_0-\delta)];\mathcal{H})$. This completes the proof by induction.
\qed

We have shown the existence and uniqueness of the mild solution $X\in\mathcal{C}^0(\mathbb{R}_+;\mathcal{H})$ for the closed-loop system associated with any initial condition $X_0 \in \mathcal{H}$ and any delay $D \in \mathcal{C}^0(\mathbb{R}_+\times[0,1];\mathbb{R})$ with $\vert D - D_0 \vert \leq \delta < D_0$. Moreover, as $v \in \mathcal{C}^0(\mathbb{R}_+;\mathcal{H})$, then the spectral reduction reported in Section~\ref{sec: problem setting} holds true. Now, the proof of the stability result stated in Theorem~\ref{thm: main result} is completed in three steps. First, a small gain argument is used to assess the stability of the truncated model (\ref{eq: spectral reduction - vect ODE}). Second, the stability of the residual infinite-dimensional dynamics (\ref{eq: spectral reduction - scalar ODE - 2}) is investigated. Finally, we will be in position to prove the stability of the closed-loop infinite-dimensional system.

\subsection{Stability analysis of the closed-loop truncated model}

The stability analysis takes the form of a small gain argument. This approach is inspired by the seminal work~\cite{karafyllis2013delay} dealing with the robustness of constant-delay predictor feedback w.r.t. uncertain and time-varying input delays.

\textbf{Step 1: use of the Artstein transformation.}
We first introduce the change of variable~\cite{artstein1982linear}:
\begin{equation}\label{eq: Artstein tranformation}
z(t) = x(t) + \int_{t-D_0}^{t} e^{(t-D_0-s)\Lambda} w(s) \,\mathrm{d}s .
\end{equation}
In particular we have from (\ref{eq: control input w}) that $w = \varphi K z$ with $z \in \mathcal{C}^1(\mathbb{R}_+;\mathbb{R}^N)$ satisfying
\begin{equation}\label{eq: ODE z - control input w}
\dot{z}(t) = \Lambda z(t) + e^{-D_0 \Lambda} w(t) + \Delta(t)
\end{equation}
for all $t \geq 0$, and thus
\begin{equation}\label{eq: ODE z - t geq t0}
\dot{z}(t) = A_\mathrm{cl} z(t) + \Delta(t)
\end{equation}
for all $t \geq t_0$.

\textbf{Step 2: estimation of $\sup_{s \in [t_0+D_0+\delta,t]} e^{\kappa s} \Vert \Delta(s) \Vert$.}
We infer that, for all $t \geq 0$,
\begin{align}
& \vert \Delta_n(t) \vert \nonumber \\
& \leq \sum\limits_{k=1}^{N} \int_0^1 \rho(\xi) \vert w_k(t-D(t,\xi)) - w_k(t-D_0) \vert \vert e_k(\xi) \vert \vert e_n(\xi) \vert \,\mathrm{d}\xi \nonumber \\
& \leq \sum\limits_{k=1}^{N} \sup\limits_{\tau \in [D_0-\delta,D_0+\delta]} \vert w_k(t-\tau) - w_k(t-D_0) \vert \label{eq: prel estimate Delta_n}
\end{align}
where it has been used that, by Cauchy-Schwarz inequality, $\int_0^1 \rho(\xi)  \vert e_k(\xi) \vert \vert e_n(\xi) \vert \,\mathrm{d}\xi \leq \Vert e_k \Vert_\mathcal{H} \Vert e_n \Vert_\mathcal{H} = 1$. Now, as $w_k = \varphi K_k z$ with $\varphi(s) = 1$ for $s \geq t_0$, we have for all $t \geq t_0 + D_0 + \delta$ that
\begin{align*}
\vert \Delta_n(t) \vert
& \leq \sum\limits_{k=1}^{N} \Vert K_k \Vert \sup\limits_{\tau \in [D_0-\delta,D_0+\delta]} \Vert z(t-\tau) - z(t-D_0) \Vert
\end{align*}
hence
\begin{equation}\label{eq: estimate Delta - 1}
\Vert \Delta(t) \Vert
\leq C_0 \sup\limits_{\tau \in [D_0-\delta,D_0+\delta]} \Vert z(t-\tau) - z(t-D_0) \Vert
\end{equation}
for all $t \geq t_0 + D_0 + \delta$ with $C_0 = \sqrt{N} \sum_{k=1}^{N} \Vert K_k \Vert$.

As $A_\mathrm{cl}$ is Hurwitz, we consider constants $M \geq 1$ and $\sigma > 0 $ such that $\Vert e^{A_\mathrm{cl} t} \Vert \leq M e^{-\sigma t}$ for all $t \geq 0$. Integrating (\ref{eq: ODE z - t geq t0}), we obtain for all $t \geq t_0+D_0+\delta$ and $\tau \in [D_0-\delta,D_0+\delta]$ that
\begin{equation*}
z(t-\tau) = e^{A_\mathrm{cl} (D_0-\tau)} z(t-D_0) + \int_{t-D_0}^{t-\tau} e^{A_\mathrm{cl} (t-\tau-s)} \Delta(s) \,\mathrm{d}s
\end{equation*}
from which we obtain that
\begin{align*}
& \Vert z(t-\tau) - z(t-D_0) \Vert \\
& \quad\leq \Vert e^{A_\mathrm{cl} (D_0-\tau)} - I \Vert \Vert z(t-D_0) \Vert \\
& \quad\phantom{\leq}\, + \left\Vert \int_{t-D_0}^{t-\tau} e^{A_\mathrm{cl} (t-\tau-s)} \Delta(s) \,\mathrm{d}s \right\Vert \\
& \quad\leq ( e^{\Vert A_\mathrm{cl} \Vert \delta} - 1 ) \Vert z(t-D_0) \Vert \\
& \quad\phantom{\leq}\, + M_\delta \left\vert \int_{t-D_0}^{t-\tau} e^{-\sigma (t-\tau-s)} \Vert \Delta(s) \Vert \,\mathrm{d}s \right\vert .
\end{align*}
where $M_\delta = \max( M , e^{\Vert A_\mathrm{cl} \Vert \delta} )$. For any $\kappa \in (0,\sigma)$, to be specified later, we have
\begin{align*}
& \left\vert \int_{t-D_0}^{t-\tau} e^{-\sigma (t-\tau-s)} \Vert \Delta(s) \Vert \,\mathrm{d}s \right\vert \\
& \leq e^{-\sigma(t-\tau)} \left\vert \int_{t-D_0}^{t-\tau} e^{(\sigma-\kappa) s} \,\mathrm{d}s \right\vert \sup\limits_{s \in [t-(D_0+\delta),t-(D_0-\delta)]} e^{\kappa s} \Vert \Delta(s) \Vert .
\end{align*}
Moreover, one has
\begin{align*}
e^{-\sigma(t-\tau)} \left\vert \int_{t-D_0}^{t-\tau} e^{(\sigma-\kappa) s} \,\mathrm{d}s \right\vert
& \leq \dfrac{e^{-\kappa(t-D_0)}}{\sigma-\kappa} \left\vert e^{\kappa (\tau-D_0)} - e^{\sigma (\tau-D_0)} \right\vert \\
& \quad\leq  \dfrac{\sigma\delta e^{\sigma\delta}}{\sigma-\kappa} e^{-\kappa(t-D_0)} ,
\end{align*}
where the last estimate is a consequence of the mean value theorem. Combining the three latter estimates, we infer that, for all $t \geq t_0+D_0+\delta$,
\begin{align*}
& \sup\limits_{\tau\in[D_0-\delta,D_0+\delta]} \Vert z(t-\tau) - z(t-D_0) \Vert \\
& \leq ( e^{\Vert A_\mathrm{cl} \Vert \delta} - 1 ) \Vert z(t-D_0) \Vert \\
& \phantom{\leq}\, + \dfrac{M_\delta \sigma\delta e^{\sigma\delta}}{\sigma-\kappa} e^{-\kappa(t-D_0)} \sup\limits_{s \in [t-(D_0+\delta),t-(D_0-\delta)]} e^{\kappa s} \Vert \Delta(s) \Vert.
\end{align*}
Thus, we infer from (\ref{eq: estimate Delta - 1}) that, for all $t \geq t_0+D_0+\delta$,
\begin{align}
& \sup\limits_{s \in [t_0+D_0+\delta,t]} e^{\kappa s} \Vert \Delta(s) \Vert  \label{eq: estimate Delta - 2} \\
& \quad\leq C_1(\delta) \sup\limits_{s \in [t_0+\delta,t-D_0]} e^{\kappa s} \Vert z(s) \Vert \nonumber \\
& \quad\phantom{\leq}\, + C_2(\delta) \sup\limits_{s \in [t_0,t-(D_0-\delta)]} e^{\kappa s} \Vert \Delta(s) \Vert \nonumber
\end{align}
where $C_1(\delta) = C_0 e^{\kappa D_0} ( e^{\Vert A_\mathrm{cl} \Vert \delta} - 1 )$ and $C_2(\delta) = \frac{M_\delta C_0 \sigma\delta e^{\sigma\delta}}{\sigma-\kappa} e^{\kappa D_0}$.

\textbf{Step 3: estimation of $\sup_{s \in [t_0,t]} e^{\kappa s} \Vert z(s) \Vert$.}
We now integrate (\ref{eq: ODE z - t geq t0}) over $[t_0,t]$ for $t \geq t_0$. Recalling that $0 < \kappa < \sigma$, this yields
\begin{align*}
\Vert z(t) \Vert
& \leq M e^{-\sigma(t-t_0)} \Vert z(t_0) \Vert
+ M \int_{t_0}^{t} e^{-\sigma (t-\tau)} \Vert \Delta(\tau) \Vert \,\mathrm{d}\tau \\
& \leq M e^{-\kappa(t-t_0)} \Vert z(t_0) \Vert
+ \dfrac{M}{\sigma-\kappa} e^{-\kappa t} \sup\limits_{s \in [t_0,t]} e^{\kappa s} \Vert \Delta(s) \Vert
\end{align*}
hence
\begin{equation}\label{eq: estimate z - 1}
\sup\limits_{s \in [t_0,t]} e^{\kappa s} \Vert z(s) \Vert
\leq M_\delta e^{\kappa t_0} \Vert z(t_0) \Vert
+ \dfrac{M_\delta}{\sigma-\kappa} \sup\limits_{s \in [t_0,t]} e^{\kappa s} \Vert \Delta(s) \Vert
\end{equation}
for all $t \geq t_0$.

\textbf{Step 4: exponential stability of $z(t)$.}
Combining estimates (\ref{eq: estimate Delta - 2}-\ref{eq: estimate z - 1}), we deduce that
\begin{align}
& \sup\limits_{s \in [t_0+D_0+\delta,t]} e^{\kappa s} \Vert \Delta(s) \Vert \label{eq: estimate Delta - 3} \\
& \quad \leq C_3(\delta) \Vert z(t_0) \Vert
+ \eta(\delta) \sup\limits_{s \in [t_0,t]} e^{\kappa s} \Vert \Delta(s) \Vert \nonumber
\end{align}
for all $t \geq t_0+D_0+\delta$ with $C_3(\delta) = M_\delta e^{\kappa t_0} C_1(\delta)$ and
\begin{align*}
\eta(\delta)
& = \dfrac{M_\delta C_1(\delta)}{\sigma-\kappa} + C_2(\delta) \\
& = \dfrac{M_\delta C_0 e^{\kappa D_0}}{\sigma-\kappa} \left\{ (e^{\Vert A_\mathrm{cl} \Vert \delta} - 1) + \sigma\delta e^{\sigma\delta} \right\} .
\end{align*}
From the small gain assumption (\ref{eq: small gain condition}), a continuity argument at $\kappa = 0$ shows the existence of\footnote{We recall that $\gamma > 0$ has been selected such that $\lambda_n \leq - \gamma$ for all $n \geq N+1$.} $\kappa \in (0,\min(\sigma , \gamma/2))$ such that $0 \leq \eta(\delta) < 1$. We fix such a $\kappa \in (0,\min(\sigma , \gamma/2))$ for the remainder of the proof. Noting that the supremums appearing in (\ref{eq: estimate Delta - 3}) are finite, we infer from this estimate that, for all $t \geq t_0+D_0+\delta$,
\begin{align*}
& \sup\limits_{s \in [t_0+D_0+\delta,t]} e^{\kappa s} \Vert \Delta(s) \Vert \\
& \quad \leq \dfrac{C_3(\delta)}{1-\eta(\delta)} \Vert z(t_0) \Vert
+ \dfrac{\eta(\delta)}{1-\eta(\delta)} \sup\limits_{s \in [t_0,t_0+D_0+\delta]} e^{\kappa s} \Vert \Delta(s) \Vert
\end{align*}
From (\ref{eq: estimate z - 1}) and using the estimate $\sup_{s \in [t_0,t]} e^{\kappa s} \Vert \Delta(s) \Vert \leq \sup_{s \in [t_0,t_0+D_0+\delta]} e^{\kappa s} \Vert \Delta(s) \Vert + \sup_{s \in [t_0+D_0+\delta,t]} e^{\kappa s} \Vert \Delta(s) \Vert$, we have for all $t \geq t_0+D_0+\delta$
\begin{align*}
& \sup\limits_{s \in [t_0,t]} e^{\kappa s} \Vert z(s) \Vert \\
& \quad \leq C_4(\delta) \Vert z(t_0) \Vert
+ C_5(\delta) \sup\limits_{s \in [t_0,t_0+D_0+\delta]} e^{\kappa s} \Vert \Delta(s) \Vert
\end{align*}
with  $C_4(\delta) = M_\delta \left\{ e^{\kappa t_0} + \frac{C_3(\delta)}{(\sigma-\kappa)(1-\eta(\delta))} \right\}$ and $C_5(\delta) = \frac{M_\delta}{\sigma-\kappa} \left\{ 1 + \frac{\eta(\delta)}{1-\eta(\delta)} \right\} = \frac{M_\delta}{(\sigma-\kappa)(1-\eta(\delta))}$. This yields, for all $t \geq t_0+D_0+\delta$,
\begin{align}
\Vert z(t) \Vert
& \leq C_4(\delta) e^{-\kappa t} \Vert z(t_0) \Vert \label{eq: estimate z - 2} \\
& \phantom{\leq}\, + C_5(\delta) e^{-\kappa t} \sup\limits_{s \in [t_0,t_0+D_0+\delta]} e^{\kappa s} \Vert \Delta(s) \Vert . \nonumber
\end{align}
We now evaluate, in function of the initial condition $z(0) = x(0)$, the two terms on the right hand side of (\ref{eq: estimate z - 2}). We recall that $w = \varphi K z$. On one hand we have for $t \leq D_0-\delta$ that $t - D(t,\xi) \leq t-(D_0-\delta) \leq 0$ and $t - D_0 \leq 0$. From (\ref{eq: structure control input v and v0}), we obtain $v(t)=v_0(t)=0$ and $\Delta_n(t) = \left< v(t) - v_{0}(t) , e_n \right> = 0$, hence $\Delta(t)=0$. On the other hand, we have from (\ref{eq: prel estimate Delta_n}) that, for $t > D_0-\delta$,
\begin{align*}
\vert \Delta_n(t) \vert
& \leq 2 \sum\limits_{k=1}^{N} \sup\limits_{s \in [t-(D_0+\delta),t-(D_0-\delta)]} \vert w_k(s) \vert \\
& \leq 2 \sum\limits_{k=1}^{N} \Vert K_k \Vert \sup\limits_{s \in [0,\max(t-(D_0-\delta),0)]} \Vert z(s) \Vert ,
\end{align*}
where we have used that $w_k(s) = 0$ for $s \leq 0$ and $\vert w_k(s) \vert \leq \Vert K_k \Vert \Vert z(s) \Vert$ for $s \geq 0$. In both cases, we obtain that, for all $t \geq 0$,
\begin{equation}\label{eq: rough estimate delta - 1}
\Vert \Delta(t) \Vert \leq 2 C_0 \sup\limits_{s \in [0,\max(t-(D_0-\delta),0)]} \Vert z(s) \Vert .
\end{equation}
We now show by induction that, for any $k \geq 1$, there exists $\alpha_k \geq 0$ such that $\Vert z(t) \Vert \leq \alpha_k \Vert x(0) \Vert$ for all $0 \leq t \leq k(D_0-\delta)$.

\emph{Initialization.} For $0 \leq t \leq D_0-\delta$, we have $\dot{z}(t) = \Lambda z(t) + e^{-D_0 \Lambda} w(t) = ( \Lambda + \varphi(t) e^{-D_0 \Lambda} K ) z(t)$ hence $\Vert \dot{z}(t) \Vert \leq C_6 \Vert z(t) \Vert$ with $C_6 = \Vert \Lambda \Vert + \Vert e^{-D_0 \Lambda} K \Vert$. In particular we have $\Vert z(t) \Vert \leq \Vert x(0) \Vert + C_6 \int_0^t \Vert z(s) \Vert \,\mathrm{d}s$. The application of Gr{\"o}nwall's inequality~\cite[Lem.~A.6.7]{Curtain2012} yields $\Vert z(t) \Vert \leq \alpha_1 \Vert x(0) \Vert$ for all $0 \leq t \leq D_0-\delta$ with $\alpha_1 = 1 + e^{C_6 (D_0 - \delta)}$.

\emph{Induction.} Assume that $\Vert z(t) \Vert \leq \alpha_k \Vert x(0) \Vert$ for all $0 \leq t \leq k(D_0-\delta)$. Recalling that $\dot{z}(t) = \Lambda z(t) + e^{-D_0 \Lambda} w(t) + \Delta(t) = ( \Lambda + \varphi(t) e^{-D_0 \Lambda} K ) z(t) + \Delta(t)$, we obtain from (\ref{eq: rough estimate delta - 1}) that, for all $0 \leq t \leq (k+1)(D_0-\delta)$,
\begin{align*}
\Vert \dot{z}(t) \Vert
& \leq C_6 \Vert z(t) \Vert + 2 C_0 \sup\limits_{s \in [0,\max(t-(D_0-\delta),0)]} \Vert z(s) \Vert \\
& \leq C_6 \Vert z(t) \Vert + 2 C_0 \alpha_k \Vert x(0) \Vert .
\end{align*}
The use of  Gr{\"o}nwall's inequality shows the existence of $\alpha_{k+1} \geq 1$ such that $\Vert z(t) \Vert \leq \alpha_{k+1} \Vert x(0) \Vert$ for all $0 \leq t \leq (k+1)(D_0-\delta)$. This completes the proof by induction.

Consequently, we have the existence of a constant $\alpha \geq 1$, independent of $X_0$, such that $\Vert z(t) \Vert \leq \alpha \Vert x(0) \Vert$ for all $0 \leq t \leq t_0+D_0+\delta$. Moreover, we obtain from (\ref{eq: rough estimate delta - 1}) that $\Vert \Delta(t) \Vert \leq 2C_0 \alpha \Vert x(0) \Vert$ for all $0 \leq t \leq t_0+D_0+\delta$.

We can now conclude on the exponential stability of $z$. On one hand, we have for all $0 \leq t \leq t_0+D_0+\delta$ that
\begin{equation*}
\Vert z(t) \Vert
\leq \alpha \Vert x(0) \Vert
\leq \alpha e^{\kappa(t_0+D_0+\delta)} e^{-\kappa t} \Vert x(0) \Vert .
\end{equation*}
On the other hand, we obtain from (\ref{eq: estimate z - 2}) that, for all $t \geq t_0+D_0+\delta$,
\begin{equation*}
\Vert z(t) \Vert
\leq \alpha \{ C_4(\delta) + 2 C_0 C_5(\delta) e^{\kappa(t_0+D_0+\delta)}  \} e^{-\kappa t} \Vert x(0) \Vert .
\end{equation*}
Combining the two latter estimates, we obtain the existence of a constant $C_7 \geq 0$ such that $\Vert z(t) \Vert \leq C_7 e^{-\kappa t} \Vert x(0) \Vert$ for all $t \geq 0$.

\textbf{Step 5: exponential stability of the system in its original coordinates.}
Recalling that $w = \varphi K z$ with $0 \leq \varphi \leq 1$, we infer that
\begin{equation}\label{eq: final estimate u}
\Vert u(t) \Vert_\mathcal{H} = \Vert w(t) \Vert
\leq C_7 \Vert K \Vert e^{-\kappa t} \Vert x(0) \Vert
\end{equation}
for all $t \geq 0$. Finally, we obtain from (\ref{eq: Artstein tranformation}) that, for all $t \geq 0$,
\begin{align}
\Vert x(t) \Vert
& \leq \Vert z(t) \Vert
+ \int_{t-D_0}^{t} e^{\vert t-D_0-s \vert \Vert \Lambda \Vert} \Vert w(s) \Vert \,\mathrm{d}s \label{eq: final estimate x} \\
& \leq C_8 e^{-\kappa t} \Vert x(0) \Vert \nonumber
\end{align}
with $C_8 = ( 1 + \Vert K \Vert e^{(\kappa+\Vert \Lambda \Vert)D_0}/\kappa ) C_7$. Thus we have shown the exponential stability of the system trajectories, as well as the exponential decay of the control input, for the closed-loop truncated model (\ref{eq: spectral reduction - vect ODE}).

\subsection{Stability analysis of the residual infinite-dimensional dynamics}\label{subsec: stab residual dynamics}

We now investigate the stability of the residual infinite-dimensional dynamics (\ref{eq: spectral reduction - scalar ODE - 2}). We consider in this subsection integers $n \geq N+1$ for which we recall that $\lambda_n \leq - \gamma < 0$. We also recall that $\kappa >0 $ has been selected such that $0 < 2\kappa < \gamma$. Now, integrating (\ref{eq: spectral reduction - scalar ODE - 2}), we infer that $x_n(t) = e^{\lambda_n t} x_n(0) + \int_0^t e^{\lambda_n (t-\tau)} v_n(\tau) \,\mathrm{d}\tau$ for all $t \geq 0$. Thus we have
\begin{align*}
\vert x_n(t) \vert^2
& \leq 2 e^{-2\gamma t} \vert x_n(0) \vert^2 + 2 \left\{ \int_0^t e^{-\gamma(t-\tau)} \vert v_n(\tau) \vert \,\mathrm{d}\tau \right\}^2 \\
& \leq 2 e^{-2\gamma t} \vert x_n(0) \vert^2 + \dfrac{2}{\gamma} \int_0^t e^{-\gamma(t-\tau)} \vert v_n(\tau) \vert^2 \,\mathrm{d}\tau
\end{align*}
where the latter estimate is obtained by using Cauchy-Schwarz inequality. Summing the latter estimate for $n \geq N+1$, we deduce that
\begin{align}
\sum\limits_{n \geq N+1} \vert x_n(t) \vert^2
& \leq 2 e^{-2\gamma t} \sum\limits_{n \geq N+1} \vert x_n(0) \vert^2 \label{eq: prel estimate inf dim dynamics} \\
& \phantom{\leq}\, + \dfrac{2}{\gamma} \int_0^t e^{-\gamma(t-\tau)} \Vert v(\tau) \Vert_\mathcal{H}^2 \,\mathrm{d}\tau \nonumber
\end{align}
for all $t \geq 0$. We now need to evaluate the term $\Vert v(\tau) \Vert_\mathcal{H}$. From (\ref{eq: structure control input v}), we have
\begin{equation*}
\Vert v(t) \Vert_\mathcal{H} \leq \sum\limits_{k=1}^{N} \sqrt{\int_0^1 \rho(\xi) \vert w_k(t-D(t,\xi)) e_k(\xi) \vert^2 \,\mathrm{d}\xi} .
\end{equation*}
Noting from (\ref{eq: final estimate u}) that
\begin{align*}
\vert w_k(t-D(t,\xi)) \vert
& \leq \Vert w(t-D(t,\xi)) \Vert \\
& \leq \sup\limits_{\tau\in[t-(D_0+\delta),t-(D_0-\delta)]} \Vert w(\tau) \Vert \\
& \leq C_7 \Vert K \Vert e^{\kappa (D_0+\delta)} e^{-\kappa t} \Vert x(0) \Vert ,
\end{align*}
we obtain from the two latter estimates that, for all $t \geq 0$,
\begin{equation*}
\Vert v(t) \Vert_\mathcal{H} \leq N C_7 \Vert K \Vert e^{\kappa (D_0+\delta)} e^{-\kappa t} \Vert x(0) \Vert ,
\end{equation*}
were we have used that $e_k \in \mathcal{H}$ is a unit vector. Since $0 < 2 \kappa < \gamma$, we have the following estimate:
\begin{align*}
\int_0^t e^{-\gamma(t-\tau)} e^{-2\kappa \tau} \,\mathrm{d}\tau
& = e^{-\gamma t} \int_0^t e^{(\gamma-2\kappa) \tau} \,\mathrm{d}\tau
\leq \dfrac{1}{\gamma-2\kappa} e^{-2\kappa t} .
\end{align*}
Using the two latter estimates into (\ref{eq: prel estimate inf dim dynamics}), we infer that
\begin{align}
\sum\limits_{n \geq N+1} \vert x_n(t) \vert^2
& \leq 2 e^{-2\kappa t} \sum\limits_{n \geq N+1} \vert x_n(0) \vert^2 \label{eq: estimate inf dim dynamics} \\
& \phantom{\leq}\, + C_9^2 e^{-2\kappa t} \Vert x(0) \Vert^2 \nonumber
\end{align}
for all $t \geq 0$, where $C_9 \geq 0$ is given by $C_9^2 = \frac{2 N^2 C_7^2 \Vert K \Vert^2 e^{2\kappa (D_0+\delta)}}{\gamma(\gamma-2\kappa)}$.

\subsection{Conclusion of the proof of the main result}

Combining estimates (\ref{eq: final estimate x}) and (\ref{eq: estimate inf dim dynamics}), we thus infer that, for all $t \geq 0$,
\begin{align*}
\Vert X(t) \Vert_\mathcal{H}^2
& = \Vert x(t) \Vert^2 + \sum\limits_{n \geq N+1} \vert x_n(t) \vert^2
\leq C_{10}^2 e^{-2\kappa t} \Vert X_0 \Vert_\mathcal{H}^2
\end{align*}
where $C_{10} \geq 0$ is given by $C_{10}^2 = \max (2,C_8^2+C_9^2)$. Recalling that the command input $u$ satisfies the estimate (\ref{eq: final estimate u}), this completes the proof of the main result.

\section{Numerical example}\label{sec: numerical application}

We illustrate the result of Theorem~\ref{thm: main result} based on the reaction-diffusion system described by (\ref{eq: pb setting}) in the case $\rho = 1$, $p = 0.015$, $q = 0.35$, $\theta_1 = \pi/3$, and $\theta_2 = \pi/10$. The open-loop system is unstable with $\lambda_1 \approx 0.317$ and $\lambda_2 \approx 0.116$ while all other modes are stable with $\lambda_3 \approx -0.342$. Thus we set $N=2$.  We consider the nominal value of the delay $D_0 = 1\,\mathrm{s}$. We impose the location $-0.3$ for the two poles of the closed-loop truncated dynamics. In this case, the small gain condition (\ref{eq: small gain condition}) is satisfied for $\delta = 0.237$, allowing to apply the stability result stated in Theorem~\ref{thm: main result}.

\begin{figure}
     \centering
		\includegraphics[width=3in]{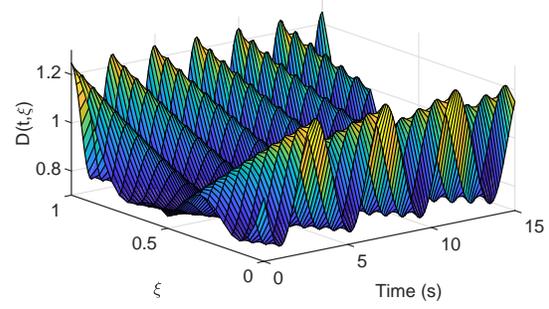}
		\label{fig: sim - delay}
     \caption{Time and spatial evolution of the input delay $D(t,\xi)$}
     \label{fig: delay}
\end{figure}

\begin{figure}
     \centering
     	\subfigure[State $y(t,\xi)$]{
		\includegraphics[width=3in]{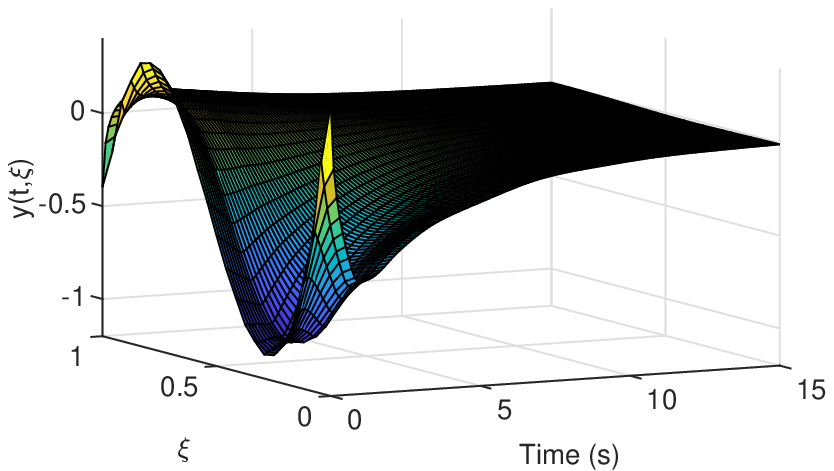}
		\label{fig: sim - state}
		}
     	\subfigure[Delayed control $v(t,\xi) = u(t-D(t,\xi),\xi)$]{
		\includegraphics[width=3in]{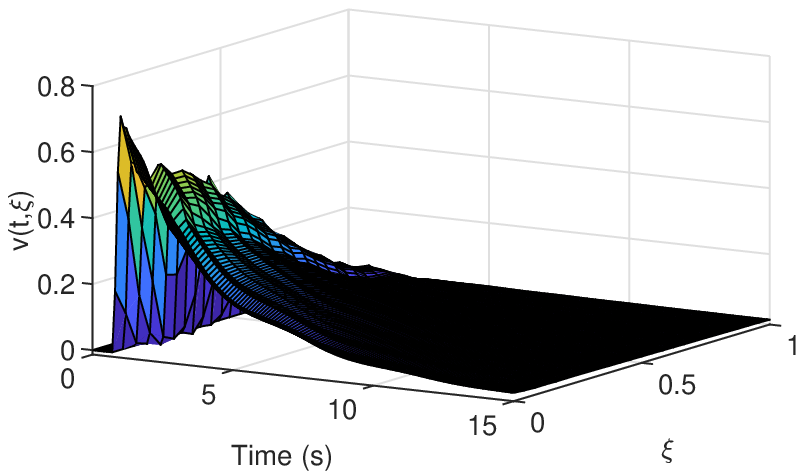}
		\label{fig: sim - distributed input}
		}
     	\subfigure[Control input $w(t)$]{
		\includegraphics[width=3in]{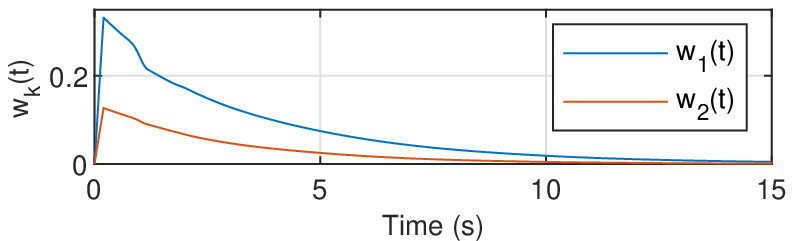}
		\label{fig: sim - inputs}
		}
     \caption{Time evolution of the closed-loop system}
     \label{fig: sim}
\end{figure}

For simulation purposes, we consider the time- and spatially-varying distributed input delay $D(t,\xi) = 0.77 + 0.23 \vert 2\xi-1 \vert \left\{ 1 + \sin([3/2+\xi]t+[11\xi-3]) \right\}$ for $t \geq 0$ and $\xi \in [0,1]$; see Fig.~\ref{fig: delay}. In particular, we have that $\vert D - D_0 \vert \leq 0.23 \leq \delta$. The initial condition is selected as $y_0(\xi) = (1-2\xi)/2 + 20 \xi (1-\xi) (\xi-3/5)$. We set the transition time as $t_0 = 0.2\,\mathrm{s}$ with $\varphi$ linearly increasing from $0$ to $1$ on $[0,t_0]$. The numerical scheme consists of the modal approximation of the reaction-diffusion equation by its $20$ dominant modes. The solution of the the implicit equation (\ref{eq: control input w}), used to implement the feedback law (\ref{eq: structure control input u}), is computed based the approximation of the integral appearing in (\ref{eq: control input w}) by a Riemann sum. The corresponding simulation results are depicted in Fig~\ref{fig: sim}. They are compliant with the predictions of Theorem~\ref{thm: main result}.

\section{Conclusion}\label{sec: conclusion}
This paper discussed the problem of in-domain stabilization of a class of infinite-dimensional systems, which operate on a weighted space of square integrable functions over a compact interval, in the presence of an uncertain time- and spatially-varying delay in the distributed actuation. This class includes, for example, reaction-diffusion PDEs. The spatially-varying nature of the delay induces new challenges because it introduces a strong coupling between the space and time variables compared to only time-varying delays configurations. We solved this control design problem by synthesizing a constant-delay predictor feedback on a finite-dimensional truncated model capturing the unstable modes of the original plant. Invoking a small gain argument, we showed that the resulting closed-loop system is exponentially stable provided the fact that the deviations of the delay around its nominal value are small enough. As small gain conditions are, in general, conservative, future works will be devoted to the derivation of relaxed stability conditions.


\bibliographystyle{plain}        
\bibliography{autosam}           




\end{document}